\documentclass{article}
\usepackage[letterpaper, portrait, margin=2.5cm]{geometry}

\usepackage{verbatim}
\usepackage[utf8]{inputenc}
\usepackage{amsmath,amssymb,amsthm}
\usepackage{tikz-cd}
\usepackage{bbm}

\usepackage{comment}
\usepackage{stmaryrd}
\usepackage[new]{old-arrows}
\usepackage{graphicx}
\graphicspath{ {./Images/} }
\usepackage{todonotes}

\usepackage [english]{babel}

\usepackage{hyperref}

\newcommand{\AAA}{\mathbb A}

\newcommand{\CC}{\mathbb C}

\newcommand{\FF}{\mathbb F}

\newcommand{\PP}{\mathbb P}

\newcommand{\VV}{\mathbb V}

\newcommand{\ZZ}{\mathbb Z}

\newcommand{\mcc}{\mathcal C}

\newcommand{\lbm}{\left[ \begin{matrix}}
\newcommand{\rem}{\end{matrix} \right]}
\newcommand{\lp}{\left(}
\newcommand{\rp}{\right)}
\newcommand{\lb}{\left\{}
\newcommand{\rb}{\right\}}

\newcommand{\ra}{\rightarrow}

\DeclareMathOperator{\PSL}{PSL}

\DeclareMathOperator{\RD}{RD}
\DeclareMathOperator{\ed}{ed}

\DeclareMathOperator{\MOO}{M_{11}}
\DeclareMathOperator{\MOW}{M_{12}}
\DeclareMathOperator{\MWW}{M_{22}}
\DeclareMathOperator{\MWH}{M_{23}}
\DeclareMathOperator{\MWF}{M_{24}}
\DeclareMathOperator{\COO}{Co_1}
\DeclareMathOperator{\COW}{Co_2}
\DeclareMathOperator{\COH}{Co_3}
\DeclareMathOperator{\JO}{J_1}
\DeclareMathOperator{\JW}{J_2}
\DeclareMathOperator{\JH}{J_3}
\DeclareMathOperator{\JF}{J_4}
\DeclareMathOperator{\SUZ}{Suz}
\DeclareMathOperator{\MCL}{McL}
\DeclareMathOperator{\HS}{HS}
\DeclareMathOperator{\FIWW}{Fi_{22}}
\DeclareMathOperator{\FIWH}{Fi_{23}}
\DeclareMathOperator{\FIWF}{Fi_{24}'}
\DeclareMathOperator{\THO}{Th}
\DeclareMathOperator{\HN}{HN}
\DeclareMathOperator{\HE}{He}
\DeclareMathOperator{\B}{B}
\DeclareMathOperator{\MO}{M}
\DeclareMathOperator{\ON}{O'N}
\DeclareMathOperator{\RU}{Ru}
\DeclareMathOperator{\LY}{Ly}

\theoremstyle{definition}
\newtheorem{definition}{Definition}[section]

\newtheorem{question}[definition]{Question}

\title{A Summary of Known Bounds on the Essential Dimension and Resolvent Degree of Finite Groups}
\author{Alexander J. Sutherland}
\date{As of \today}

\begin{document}

\maketitle

\begin{abstract}
We summarize what is currently known about $\ed(G)$ and $\RD(G)$ for finite groups $G$ over $\CC$ (i.e. in characteristic 0). In Appendix \ref{app:The Case of PSL(2,11)}, we also give an argument which improves the known bound on $\RD(\PSL(2,\FF_{11}))$.
\end{abstract}

\setcounter{tocdepth}{1}
\tableofcontents

\section{Introduction}

The art of solving polynomials has a long and storied history. These investigations have led to a spectrum of modern frameworks to classify the complexity of phenomena in algebra and geometry; the ends of this spectrum are \textbf{essential dimension} and \textbf{resolvent degree}. 

Using Galois theory, we can reformulate the classical problem
\begin{center}
	\emph{``Determine a simplest formula for the general degree $n$ polynomial.''};
\end{center}

\noindent
into
\begin{center}
	\emph{``Determine $\ed(S_n)$ and $\RD(S_n)$, the essential dimension and resolvent degree of $S_n$.''};
\end{center}

\noindent
where $S_n$ denotes the symmetric group on $n$ letters. 

Essential dimension and resolvent degree are not defined only for the symmetric group, but for all finite groups. Note that for any finite group $G$, \cite[Lemma 3.2]{FarbWolfson2019} yields
\begin{equation*}
	\RD(G) \leq \ed(G) < \infty. 
\end{equation*}

\noindent
Consequently, we are naturally led to the following question.

\begin{question}[Main Question]\label{ques:Main Question}
For which finite groups does $\RD(G) = \ed(G)$ and for which finite groups is $\RD(G) < \ed(G)$?
\end{question}

In this short note, we will address what is currently known about Question \ref{ques:Main Question}. In Section \ref{sec:Upper Bounds on RD(G)}, we discuss the finite groups $G$ with explicit upper bounds on $\RD(G)$ in the literature. In Section \ref{sec:Lower Bounds on ed(G)}, we discuss what is known about $\ed(G)$ for the groups given in Section 1. In Section \ref{sec:Summary and Tables}, we summarize the answers to Question \ref{ques:Main Question}. Note that we will not give proofs of the results, but we will indicate where these results can be found.

\paragraph{Additional Context} 

For those interested in the history leading to resolvent degree, we refer the reader to \cite[Chapter 1]{Sutherland2022}. For those interested in a larger-scale perspective on essential dimension and resolvent degree, we refer the reader to the survey of Reichstein \cite{Reichstein2021}.

For the expert reader, we note that essential dimension and resolvent degree can be (and often are) considered over arbitrary base fields and for arbitrary algebraic groups; see \cite{EdensReichstein2023}, for example, where the essential dimension of the symmetric groups is considered in prime characteristic. However, we restrict our attention here to the classical case where our base field is $\CC$ and $G$ is a finite group.

\paragraph{Acknowledgements}

I would like to thank Daniel Litt for recently asking what is known about Question \ref{ques:Main Question}, which prompted this expository note. I would like to thank Jesse Wolfson for several helpful comments on a draft, including providing information on several relevant references.

\section{Upper Bounds on $\RD(G)$}\label{sec:Upper Bounds on RD(G)}

For a finite group $G$ with simple factors $G_1,\dotsc,G_s$, \cite[Theorem 3.3]{FarbWolfson2019} yields that
\begin{equation*}
	\RD(G) \leq \max\lb \RD(G_1),\dotsc,\RD(G_s) \rb.
\end{equation*}

Consequently, the literature has focused on $\RD(G)$ when $G$ is a finite simple group. More specifically, there are non-trivial upper bounds on $\RD(G)$ in the literature for a finite simple group $G$ in the following cases:
\begin{enumerate}
\item $G$ is a cyclic group of prime order;
\item $G$ is an alternating group $(A_n, n \geq 5)$;
\item $G$ is a simple factor of a Weyl group of type $E_6$, $E_7$, or $E_8$;
\item $G = \PSL(2,\FF_7)$ or $\PSL(2,\FF_{11})$; or
\item $G$ is a sporadic group.
\end{enumerate}

\noindent
Let us now consider each case individually.

\paragraph{Case 1:} It is immediate that $\RD(G)=1$; it is also a special case of \cite[Corollary 3.4]{FarbWolfson2019}. Indeed, $\RD(A)=1$ for any finite abelian group $A$.

\paragraph{Case 2:} Before stating known bounds, we note that $\RD(S_n) = \RD(A_n)$, whereas $\ed\lp A_n \rp \not= \ed(S_n)$; we will say more in Section \ref{sec:Lower Bounds on ed(G)}. Bring showed that $\RD(A_n) \leq n-4$ for $n \geq 5$ in \cite{Bring1786}. Segre gave a complete proof that $\RD(A_n) \leq n-5$ for $n \geq 9$ in \cite{Segre1945}, building upon the work of Hilbert \cite{Hilbert1927}. Indeed, all of our bounds come in the form $\RD(A_n) \leq n-m$ for $n$ greater than some threshold determined by $m$. For most values of $m \geq 6$, the bounds are given by $G(m)$ in \cite{Sutherland2021}; see Theorems 3.7, 3.10, and 3.27, or Appendices 5.1 and 5.2 for a summary. For $m \in [13,17] \cup [22,25]$, improvements to $G(m)$ are established in \cite[Theorem 1.1]{HeberleSutherland2023}.

\paragraph{Case 3:} We denote the simple factors of $W(E_6)$, $W(E_7)$, and $W(E_8)$ by $W(E_6)^+, W(E_7)^+, W(E_8)^+$, respectively, and note that in each case $\RD\lp W(E_n) \rp = \RD\lp W(E_n)^+ \rp$. In \cite{FarbWolfson2019}, Theorem 8.2 establishes $\RD\lp W(E_6)^+ \rp \leq 3$. In \cite{Reichstein2022}, Proposition 15.1 establishes $\RD\lp W(E_7)^+ \rp \leq 4$ and $\RD\lp W(E_8)^+ \rp \leq 5$. 

\paragraph{Case 4:} The group $\PSL(2,\FF_7)$ was investigated classically by Fricke, Klein, and Gordan. We refer the reader to \cite[Proposition 4.2.4]{FarbKisinWolfson2022}, which establishes that $\RD(\PSL(2,\FF_7)) = 1$. Using modern language, \cite{Klein1879} establishes that $\RD(\PSL(2,\FF_{11})) \leq 3$. We will say more in Appendix \ref{app:The Case of PSL(2,11)}, where we outline a proof that $\RD(\PSL(2,\FF_{11})) \leq 2$.

\paragraph{Case 5:}  The resolvent degree of the sporadic groups was recently investigated by the author in joint work with G\'{o}mez-Gonz\'{a}les and Wolfson. Corollary 4.9 of \cite{GomezGonzalesSutherlandWolfson2023} gives the following bounds:

\begin{align*}
	\RD(\JW) &\leq 5, &\RD(\MWF) &\leq 18, &\RD(\HE) &\leq 48, &\RD(\FIWH) &\leq 776,\\
	\RD(\MOO) &\leq 6, &\RD(\HS) &\leq 18, &\RD(\JO) &\leq 51, &\RD(\FIWF) &\leq 779,\\
	\RD(\MOW) &\leq 7, &\RD(\MCL) &\leq 19, &\RD(\FIWW) &\leq 74, &\RD(\JF) &\leq 1328,\\
	\RD(\MWW) &\leq 8, &\RD(\COH) &\leq 20, &\RD(\HN) &\leq 129, &\RD(\LY) &\leq 2475,\\
	\RD(\SUZ) &\leq 10, &\RD(\COW) &\leq 20, &\RD(\THO) &\leq 244, &\RD(\B) &\leq 4365,\\
	\RD(\JH) &\leq 16, &\RD(\COO) &\leq 21, &\RD(\ON) &\leq 338, &\RD(\MO) &\leq 196874.\\
	\RD(\MWH) &\leq 17, &\RD(\RU) &\leq 26.
\end{align*}

\section{Lower Bounds on $\ed(G)$}\label{sec:Lower Bounds on ed(G)}

Let us again address the five cases outlined in Section \ref{sec:Upper Bounds on RD(G)}, where non-trivial upper bounds on $\RD(G)$ are known. The author is not aware of any lower bounds on the essential dimension of the sporadic groups, so we omit Case 5 here.

\paragraph{Case 1:} We follow \cite[Section 6.1]{BuhlerReichstein1997} and recall that the rank of a finite abelian group $A$ is the minimal number of elements that generate $A$; equivalently, it is the largest $r$ such that $\lp \ZZ/p\ZZ \rp^r$ embeds into $G$ for some prime $p$. Theorem 6.1 of \cite{BuhlerReichstein1997} then yields that
\begin{equation*}
	\ed(A) = \text{rank}(A).
\end{equation*}

\noindent
Consequently, when $G \cong \ZZ/p\ZZ$, we see that $\ed(G)=1$.

\paragraph{Case 2:} We will discuss both the alternating groups $A_n$ and the symmetric groups $S_n$. It was known classically, e.g. by Felix Klein, that
\begin{align*}
 \ed(A_2) = \ed(A_3) &= 1, &\ed(A_4) &= \ed(A_5) = 2, &\ed(A_6) &= 3.
\end{align*}

\noindent
In \cite{BuhlerReichstein1997}, Theorem 6.7 further establishes that
\begin{enumerate}
	\item $\ed(A_{n+4}) \geq \ed(A_n)+2$ for $n \geq 4$,
	\item $\ed(A_n) \geq 2\lfloor \frac{n}{4} \rfloor$ for $n \geq 4$. 
\end{enumerate}

We now turn to the case of the symmetric groups. Again, it was known classically that
\begin{align*}
 \ed(S_2) = \ed(S_3) &= 2, &\ed(S_4) &= \ed(S_5) = 2.
\end{align*}

\noindent
Theorem 6.5 of \cite{BuhlerReichstein1997} additionally establishes
\begin{enumerate}
	\item $\ed(S_6)=3$,
	\item $\ed(S_n) \leq n-3$ for any $n \geq 5$,
	\item $\ed(S_{n+2}) \geq \ed(S_n)+1$ for any $n \geq 1$,
	\item $\ed(S_n) \geq \lfloor \frac{n}{2} \rfloor$ for any $n \geq 1$.
\end{enumerate}

\noindent
In \cite{Duncan2010}, Duncan showed that $\ed(A_7) = \ed(S_7) = 4$ and thus
\begin{equation*}
	\ed(S_n) \geq \left\lfloor \frac{n+1}{2} \right\rfloor, \text{ for } n \geq 7. 
\end{equation*}

As is noted in \cite[Section 1]{EdensReichstein2023}, ``the exact value of $\ed(S_n)$ is open for every $n \geq 8$, though it is widely believed that $\ed(S_n)$ should be $n-3$ for every $n \geq 5$.''

\paragraph{Case 3:} Duncan classifies the finite groups of essential dimension 2 in \cite[Theorem 1.1]{Duncan2013} and, notably, $\ed(\PSL(2,\FF_7))=2$. Proposition 4 of \cite{Beauville2014} gives the following list of possible finite simple groups of essential dimension 3: they are $A_6$ and $\PSL\lp 2,\FF_{11} \rp$. More specifically, we know $\ed(A_6)=3$ (as covered above) and $\ed(\PSL(2,\FF_{11})) = $ 3 or 4.

\paragraph{Case 4:}  The essential dimension (and essential $p$-dimension, which we do not discuss here) of finite pseudo-relection groups is determined in in \cite[Theorem 1.3]{DuncanReichstein2014}. Notably, they establish that $\ed(W(E_6)) = 4$, $\ed(W(E_7)=7$, and $\ed(W(E_8))=8$. The first equality is given explicitly in Theorem 1.3; the second is given in Remark 7.2; and the third from Example 1.2, as $\ed(W(E_8)) = \dim(V) = a(2)$.

\section{Summary and Tables}\label{sec:Summary and Tables}
We now summarize the numerical results stated above in tables for Cases 1, 2, 3, and 4.

\paragraph{Case 1:} Given that the essential dimension of an arbitrary abelian group $A$ is the essential dimension of a elementary $p$-subgroup of maximal rank, we focus on the cases of cyclic $p$-groups and elementary $p$-groups. 

\begin{center}
\begin{tabular}{|c|c|c|}
\hline
$A$ &$\RD(A)$ &$\ed(A)$\\
\hline
$\ZZ/p\ZZ$ &1 &1\\
\hline
& & \\[-9pt]
$\lp \ZZ/p\ZZ \rp^r$ &1 &$r$\\
\hline
\end{tabular}
\end{center}

\noindent
Note that $\ed\lp (\ZZ/p\ZZ)^r \rp - \RD\lp (\ZZ/p\ZZ)^r \rp = r-1$, so we can find groups where the gap between essential dimension and resolvent degree is arbitrarily large.

\paragraph{Case 2:} We record the relevant values for $A_n$ and $S_n$, up to $n=22$.

\begin{center}
\begin{tabular}{|c|c|c|c|c|}
\hline
$n$ &$\RD(A_n)=\RD(S_n)$ &$\ed(A_n)$ &$\ed(S_n)$, known &$\ed(S_n)$, conjectured\\
\hline
2 &1 &1 &1 &N/A\\
\hline
3 &1 &1 &1 &N/A\\
\hline
4 &1 &2 &2 &N/A\\
\hline
5 &1 &2 &2 &N/A\\
\hline
6 &$\leq 2$ &3 &3 &N/A\\
\hline
7 &$\leq 3$ &4 &4 &N/A\\
\hline
8 &$\leq 4$ &$\geq 4$ &$\geq 4$ &5\\
\hline
9 &$\leq 4$ &$\geq 4$ &$\geq 4$ &6\\
\hline
10 &$\leq 5$ &$\geq 4$ &$\geq 5$ &7\\
\hline
11 &$\leq 6$ &$\geq 4$ &$\geq 5$ &8\\
\hline
12 &$\leq 7$ &$\geq 6$ &$\geq 6$ &9\\
\hline
13 &$\leq 8$ &$\geq 6$ &$\geq 6$ &10\\
\hline
14 &$\leq 9$ &$\geq 6$ &$\geq 7$ &11\\
\hline
15 &$\leq 10$ &$\geq 6$ &$\geq 7$ &12\\
\hline
16 &$\leq 11$ &$\geq 8$ &$\geq 8$ &13\\
\hline
17 &$\leq 12$ &$\geq 8$ &$\geq 8$ &14\\
\hline
18 &$\leq 13$ &$\geq 8$ &$\geq 9$ &15\\
\hline
19 &$\leq 14$ &$\geq 8$ &$\geq 9$ &16\\
\hline
20 &$\leq 15$ &$\geq 10$ &$\geq 10$ &17\\
\hline
21 &$\leq 15$ &$\geq 10$ &$\geq 10$ &18\\
\hline
22 &$\leq 16$ &$\geq 10$ &$\geq 11$ &19\\
\hline
\end{tabular}
\end{center}

\noindent
If true, the conjecture that $\ed(S_n) = n-3$ for $n \geq 5$ would yield that
\begin{itemize}
\item $\ed(S_n) - \RD(S_n) \geq 1$ for $n \geq 5$,
\item $\ed(S_n) - \RD(S_n) \geq 2$ for $n \geq 9$,
\item $\ed(S_n) - \RD(S_n) \geq 3$ for $n \geq 21$,
\end{itemize}

\noindent
and so on. In fact, given this conjecture, work as old as \cite{Hamilton1836} implies that
\begin{equation*}
 \lim\limits_{n \ra \infty} \lp \ed(S_n) - \RD(S_n) \rp = \infty.
\end{equation*}

\noindent
For modern references which improve upon Hamilton's results, see \cite[Theorem 1.1]{Wolfson2021} or \cite[Theorem 1.3]{Sutherland2021}.

\paragraph{Case 3:} For $\PSL(2,\FF_7)$ and $\PSL(2,\FF_{11})$, we have the following:
\begin{center}
\begin{tabular}{|c|c|c|}
\hline
$G$ &$\RD(G)$ &$\ed(G)$\\
\hline
$\PSL(2,\FF_7)$ &1 &2\\
\hline
$\PSL(2,\FF_{11})$ &$\leq 2$ &3 or 4\\
\hline
\end{tabular}
\end{center}

\paragraph{Case 4:} For the Weyl groups $W(E_6)$, $W(E_7)$, and $W(E_8)$, we note the following bounds:
\begin{center}
\begin{tabular}{|c|c|c|}
\hline
$G$ &$\RD(G)$ &$\ed(G)$\\
\hline
& & \\[-9pt]
$W(E_6)$ &$\leq 3$ &4\\
\hline
& & \\[-9pt]
$W(E_7)$ &$\leq 4$ &7\\
\hline
& & \\[-9pt]
$W(E_8)$ &$\leq 5$ &8\\
\hline
\end{tabular}
.
\end{center}

\appendix

\section{The Case of $\PSL(2,\FF_{11})$}\label{app:The Case of PSL(2,11)}

In \cite{Klein1879}, Klein investigates $\PSL(2,\FF_{11})$. To do so, Klein considers the four-dimensional projective representation $\PP^4$ of $\PSL(2,\FF_{11})$ and constructs invariants of degrees 3 ($\nabla$), 5 ($H$), and 11 ($C$).  Using the modern computer algebra system \texttt{GAP}, one can directly compute the Molien series for the linear representation $\AAA^5$ as the following rational expression:

\begin{equation*}
	\frac{ 1 - z^2 + z^7 + z^8 - z^{11} + z^{14} + z^{15} - z^{20} + z^{22} }{ \lp 1-z^{11} \rp \lp 1-z^6 \rp \lp 1-z^5 \rp \lp 1-z^3 \rp \lp 1-z^2 \rp }.
\end{equation*}

\noindent
The series expansion of the Molien series begins:
\begin{equation*}
	1 + z^3 + z^5 + 2z^6 + z^7 + 2z^8 + 3z^9 + 3z^{10} + 4z^{11} + 6z^{12} + \cdots.
\end{equation*}

Note that $\VV(\nabla) \subseteq \PP^4$ is a cubic hypersurface in $\PP^4$ and thus we know classically that it has dense solvable points. This is enough to establish the Klein's claim that
\begin{equation*}
	\RD(\PSL(2,\FF_{11})) \leq \dim\lp \VV(\nabla) \rp = 3.
\end{equation*}

We additionally outline a quick proof of an improved bound. Before we begin, we refer the reader to \cite[Section 2]{Sutherland2021} for more on polar cones and \cite[Example 2.16]{GomezGonzalesSutherlandWolfson2023} for more on the notation $K^{(d)}$.

\begin{proof} 
For every $K$-point $P \in K(\VV(\nabla))$, the polar cone $\mcc(\VV(\nabla);P)$ has
\begin{align*}
	\deg(\mcc(\VV(\nabla);P))=6,\\
	\text{codim}(\mcc(\VV(\nabla);P)) \geq 1.
\end{align*}

\noindent
Thus, $P$ lies on a line of $\VV(\nabla)$ over $K^{(2)}$, as $\RD(A_6)=\RD(S_6) \leq 2$. For any such line $\Lambda$, we see that $\VV(\nabla) \cap \Lambda \subseteq \VV(\nabla,H)$ and $\deg\lp \VV(H) \cap \Lambda \rp = 3$. Since $\RD(A_3)=\RD(S_3)=1 \leq 2$, $\VV(\nabla,H)$ also has dense $K^{(2)}$-points and we see that
\begin{equation*}
	\RD(\PSL(2,\FF_{11})) \leq \max\lb 2,\dim(\VV(\nabla,H)) \rb = 2,
\end{equation*}

\noindent
if the action of $\PSL\lp 2,\FF_{11} \rp$ is generically free. Since $\deg\lp \VV(\nabla,H) \rp = 10$, we see that $\VV(\nabla,H)$ has at most 10 irreducible components. Since $\PSL(2,\FF_{11})$ would act on the irreducible components transitively and the smallest permuation representation of $\PSL(2,\FF_{11})$ is on 11 points, we see that $\VV(\nabla,H)$ is irreducible and thus \cite[Lemma 2.12]{GomezGonzalesSutherlandWolfson2023} yields that the action of $\PSL(2,\FF_{11})$ on $\VV(\nabla,H)$ is generically free.
\end{proof}

\end{document}